\documentclass[11pt,reqno]{amsart}
\usepackage{amsmath}
\usepackage{tikz}
\usepackage[all]{xy}
\usetikzlibrary{decorations.markings}
\usepackage{bm,mdframed}
\usepackage{color}
\usepackage{amsthm}
\usepackage{amsfonts}
\usepackage{amssymb}
\usepackage{setspace}
\usepackage[enableskew]{youngtab}
\usepackage{ytableau}
\usepackage{empheq}
\usepackage[bookmarks,colorlinks,breaklinks]{hyperref}  
\hypersetup{linkcolor=blue,citecolor=blue,filecolor=blue,urlcolor=blue} 
\theoremstyle{plain}

\definecolor{gr40}{gray}{0.40}

\newcommand{\beqn}{\begin{eqnarray}}
\newcommand{\eeqn}{\end{eqnarray}}

\newtheorem{Theorem}{Theorem}[section]
\newtheorem{Proposition}[Theorem]{Proposition}
\newtheorem{Definition}[Theorem]{Definition}
\newtheorem{Example}[Theorem]{Example}
\newtheorem{Corollary}[Theorem]{Corollary}

\newtheorem{Remark}[Theorem]{Remark}

\begin{document}
\title{Kronecker Coefficients For Some Near-Rectangular Partitions}
\author{VASU V. TEWARI}
\address{
Vasu V. Tewari\\
Department of Mathematics\\
University of British Columbia\\
Vancouver, BC V6T 1Z2\\
Canada
}
\email{\href{mailto:vasu@math.ubc.ca}{vasu@math.ubc.ca}}

\begin{abstract}
We give formulae for computing Kronecker coefficients occurring in the expansion of $s_{\mu}*s_{\nu}$, where both $\mu$ and $\nu$ are nearly rectangular, and have smallest parts equal to either 1 or 2. In particular, we study $s_{(n,n-1,1)}*s_{(n,n)}$, $s_{(n-1,n-1,1)}*s_{(n,n-1)}$, $s_{(n-1,n-1,2)}*s_{(n,n)}$, $s_{(n-1,n-1,1,1)}*s_{(n,n)}$ and $s_{(n,n,1)}*s_{(n,n,1)}$. Our approach relies on the interplay between manipulation of symmetric functions and the representation theory of the symmetric group, mainly employing the Pieri rule and a useful identity of Littlewood. As a consequence of these formulae, we also derive an expression enumerating certain standard Young tableaux of bounded height, in terms of the Motzkin and Catalan numbers.
\end{abstract}
\maketitle

An outstanding open problem in algebraic combinatorics is to derive a combinatorial formula to compute the Kronecker product of two Schur functions. Given partitions $\lambda,\mu \text{ and } \nu$, the Kronecker coefficients, $g^{\lambda}_{\mu\nu}$, occur in the decomposition of the Kronecker product $s_{\mu}*s_{\nu}$ of Schur functions in the Schur basis.
\begin{eqnarray*}
s_{\mu}*s_{\nu}=\sum_{\lambda}g^{\lambda}_{\mu\nu}s_{\lambda}
\label{kronprod1}
\end{eqnarray*} 
Alternatively, these coefficients can also be defined as the multiplicities of the irreducible representations of the symmetric group in the tensor product of two irreducible representations of the symmetric group. This interpretation immediately implies that the Kronecker coefficients are non-negative integers leading one to believe that there should be a combinatorial rule to compute these coefficients.
However, to date, there is no satisfactory positive combinatorial formula for the Kronecker product of two Schur functions. 

Besides the intrinsic interest in the problem, the motivation for discovering a combinatorial formula is the impact beyond algebraic combinatorics. For example, the Kronecker coefficients arise in quantum information theory and quantum computation \cite{ChristandlHarrowMitchison, ChristandlMitchison,W}. The problem of computing them combinatorially has got major impetus from the fact that they are of prime importance in Geometric Complexity Theory, a program of Mulmuley aimed at resolving the P vs NP problem \cite{MulmuleySohoni}. In other applications, these coefficients have been used to show the strict unimodality of $q$-binomial numbers by Pak and Panova \cite{PakPanova}. 

Attempts have been made to understand different aspects of these coefficients, for example, special cases \cite{BB,BK,BvWZ,RW,Ro,Thibon}, asymptotics \cite{BO,BO1}, stability \cite{BriandOrellanaRosas-2,V}, the complexity of computing them and conditions which guarantee that they are non-zero \cite{BI}. Recently, a combinatorial rule was given by Blasiak \cite{Blasiak} for computing $s_{\mu}*s_{\nu}$ where at least one of $\mu$ and $\nu$ is a hook shape. Finally, a certain variant, called the reduced Kronecker coefficients, has also been studied in \cite{BriandOrellanaRosas-1, BriandOrellanaRosas-2}.

The aim of this article is to derive explicit combinatorial formulae for Kronecker coefficients corresponding to partitions of \textit{near-rectangular shape}, i.e., partitions such that nearly all their parts are equal. Kronecker coefficients indexed by such partitions are conducive to manipulation, as demonstrated in \cite{BvWZ,BurgisserChristandlIkenmeyer,ChristandlHarrowMitchison, Manivel, W}. The organization of this article is as follows. In Section \ref{section:background}, we equip the reader with the required background on symmetric functions and a brief overview of relevant results. In Sections \ref{section:first case} and \ref{section: second case} we prove combinatorial formulae for the Kronecker coefficients appearing in the products $s_{(n,n-1,1)}*s_{(n,n)}$ and $s_{(n-1,n-1,1)}*s_{(n,n-1)}$ respectively. In Sections \ref{section: third case}, \ref{section: fourth case} and \ref{section: fifth case}, we state results for the products $s_{(n-1,n-1,2)}*s_{(n,n)}$, $s_{(n-1,n-1,1,1)}*s_{(n,n)}$ and $s_{(n,n,1)}*s_{(n,n,1)}$ respectively. The techniques used to obtain these results are very similar to those employed in Sections \ref{section:first case} and \ref{section: second case} and hence, the proofs are replaced by illustrative examples. The interested reader can find the complete proofs in \cite{Tewari}. Using the results obtained, we give a closed formula for the number of standard Young tableaux of height exactly $5$ and smallest part equal to $1$. This is given in Theorem \ref{theorem:main combinatorial result} in Section \ref{section: enumerative applications}. Finally, we conclude with possible future avenues in Section \ref{section:future avenues}. 

\section*{Acknowledgement}
The author would like to thank Stephanie van Willigenburg for suggesting the problem, helpful guidance and for encouraging him to write the results up.
\section{Background}\label{section:background}
We will start by defining some of the combinatorial structures that we will be encountering. All the central notions introduced in this section are covered in more detail in \cite{macdonald-1,Sa,St}. Our first definition has to do with the notion of partition.
\subsection{Partitions}
A \textit{partition} $\lambda$ is a finite list of positive integers $(\lambda_1,\ldots,\lambda_k)$ satisfying $\lambda_1\geq \lambda_2\geq \cdots\geq \lambda_k$. The integers appearing in the list are called the \textit{parts} of the partition. Given a partition $\lambda=(\lambda_1,\ldots,\lambda_k)$, the \textit{size} $\lvert \lambda\rvert$ is defined to be $\sum_{i=1}^{k}\lambda_i$. The number of parts of $\lambda$ is called the \textit{length}, and is denoted by $l(\lambda)$. If $\lambda$ is a partition satisfying $\lvert \lambda\rvert=n$, then we write it as $\lambda \vdash n$. Conventionally, there is a unique partition of size and length $0$, and we denote it by $\varnothing$.

We will be depicting a partition using its \textit{Ferrers diagram} (or \textit{Young diagram}). Given a partition $\lambda=(\lambda_1,\ldots,\lambda_k)\vdash n$, the Ferrers diagram of $\lambda$, also denoted by $\lambda$, is the left-justified array of $n$ boxes, with $\lambda_i$ boxes in the $i$-th row. We will be using the English convention, i.e. the rows are numbered from top to bottom and the columns from left to right. We refer to the box in the $i$-th row and $j$-th column by the ordered pair $(i,j)$. Finally, the \textit{transpose}, $\lambda^{t}$, of a partition $\lambda = (\lambda_{1},\lambda_{2},\ldots,\lambda_{k})$ is the partition obtained by transposing the Ferrers diagram of $\lambda$. Thus, for example, the transpose of the partition $\lambda = (5,3,3,1)$ is $\lambda^{t} = (4,3,3,1,1)$. The \textit{hooklength} associated to the box $(i,j)$, denoted by $h_{(i,j)}$ is the number $\lambda_i-j+\lambda_j^{t}-i+1$. 

If $\lambda$ and $\mu$ are partitions such that $\mu \subseteq \lambda$, i.e., $l(\mu) \leq l(\lambda)$ and  $\mu_{i} \leq \lambda_{i}$ for all $i=1,2,\ldots,l(\mu)$, then the \textit{skew shape} $\lambda / \mu$ is obtained by removing the first $\mu_{i}$ boxes from the $i$-th row of the Ferrers diagram of $\lambda$  for $1\leq i\leq l(\mu)$.  The \textit{size} of the skew shape $\lambda/\mu$, denoted by $|\lambda/\mu|$, is equal to the number of boxes in the skew shape, i.e., $\lvert \lambda \rvert-\lvert \mu \rvert$. 

Let $\mu$ and $\lambda$ be partitions. We say that $\mu\prec \lambda$ if $\mu$ can be obtained by subtracting 1 from some part of $\lambda$.
\begin{Example}
Shown below is the Ferrers diagram of $\lambda = (5,3,3,1)$ (left) and the respective hooklengths associated with each box (right).  \begin{eqnarray*}\yng(5,3,3,1) \hspace{5mm}\young(86521,532,421,1)\end{eqnarray*}
\end{Example}

Now, we will define some statistics on partitions that we will need to state our results, especially in Sections \ref{section: third case}, \ref{section: fourth case} and \ref{section: fifth case}. We will denote the number of distinct parts in a partition $\lambda$ by $d_{\lambda}$, while $d_{\lambda,2}$ will denote the number of parts of $\lambda$ from which $2$ can be subtracted so that whatever remains (once the tail of zeroes has been removed) is still partition. Finally, we denote by $R_{\lambda}$ the number of distinct parts of $\lambda$ that occur at least twice. For example, consider $\lambda=(6,5,3,3,3,2,2)$. Clearly, $d_{\lambda}=4$ and $R_{\lambda}=2$. Notice that on subtracting $2$ from the part equal to $5$ in $\lambda$, we obtain $(6,3,3,3,3,2,2)$ which is a partition while on subtracting $2$ from the rightmost part equal to $2$ in $\lambda$, we get $(6,5,3,3,3,2,0)$ which becomes a partition once we remove the 0 in the tail. Hence $d_{\lambda,2}=2$.

Given a partition $\lambda$, define a new partition $\lambda'$ as follows.
\begin{eqnarray*}
\lambda' =\left\lbrace\begin{array}{ll}\lambda & l(\lambda)\leq 4\\ (\lambda_1,\ldots,\lambda_4) & l(\lambda)>4\end{array}\right.
\end{eqnarray*}
Thus, for example, if $\lambda=(4,3,1)$ then so is $\lambda'$, but if $\lambda=(5,5,3,3,2,1)$ then $\lambda'=(5,5,3,3)$. Now given a partition $\lambda$, define $O_{\lambda}$ and $E_{\lambda}$ to be the number of odd and even parts in $\lambda'$ respectively. Also, define $O_{\lambda}^{'}$ and $E_{\lambda}^{'}$ to be the number of distinct odd parts and distinct even parts in $\lambda'$ respectively. 

To illustrate these definitions, we give an example. Consider $\lambda=(4,4,3,2,1)$. Then the number of even parts in $\lambda'=(4,4,3,2)$ is $3$. Hence $E_{\lambda}=3$, but notice the the number of distinct even parts in $\lambda'$ is just 2, i.e. $E_{\lambda}^{'}=2$. Note also that $O_{\lambda}=O_{\lambda}^{'}=1$.

\subsection{Semistandard Young tableaux}
Given partitions $\lambda$ and $\mu$ such that $\mu\subseteq \lambda$, a \textit{semistandard Young tableau} (SSYT) of \textit{shape} $\lambda/\mu$ is a filling of the boxes of the skew shape $\lambda/\mu$ with positive integers satisfying the condition that entries increase weakly along each row from left to right and increase strictly along each column from top to bottom.

A \textit{standard Young tableau} (SYT) of shape $\lambda/\mu$ is an SSYT in which the entries in the filling are distinct elements of $\{1,2,\ldots,|\lambda/\mu|\}$. We denote by $SSYT(\lambda/\mu)$ the set of all SSYT of shape $\lambda/\mu$. As a matter of convention, an SSYT of shape $\lambda/\varnothing$ will be referred to as an SSYT of shape $\lambda$. The \textit{height} of an SSYT of shape $\lambda$ is defined to be $l(\lambda)$.

The number of SYTs of shape $\lambda \vdash n$ will be denoted by $f_{\lambda}$, and it can be easily calculated by the \textit{hooklength formula} of Frame, Robinson and Thrall. 
\begin{Theorem}\cite[Theorem 1]{FRT}\label{thm: hooklength formula}
Given a partition $\lambda$ of $n$,
\begin{eqnarray*}
f_{\lambda}=\displaystyle\frac{n!}{\prod_{(i,j)\in \lambda}h_{(i,j)}}.
\end{eqnarray*}
\end{Theorem}
Finally, given a skew shape $\lambda/\mu$, we will associate a monomial $x^{T}$ to every $T \in SSYT(\lambda/\mu)$ in the following manner. 
\begin{eqnarray*}
x^{T}=\prod_{(i,j)\in \lambda/\mu}x_{T_{(i,j)}}
\end{eqnarray*}
\begin{Example}
An SSYT (left) and an SYT (right) of shape $\lambda=(4,3,1,1)$ are shown below.
\begin{eqnarray*}
\young(1333,244,5,6)\hspace{15mm}
\young(1357,248,6,9)
\end{eqnarray*}
The monomial associated with the SSYT on the left is $x_1x_2x_3^3x_4^2x_5x_6$.
\label{ssytexample}
\end{Example}

\subsection{Symmetric functions}
We will denote the algebra of symmetric functions by $\Lambda$. It is the algebra freely generated over $\mathbb{Q}$ by countably many commuting variables $\{p_1,p_2,\ldots\}$. Assigning the degree $i$ to $p_i$ (and then extending this multiplicatively) gives $\Lambda$ the structure of a graded algebra. A basis for the degree $n$ component of $\Lambda$, denoted by $\Lambda^{n}$, is given by the \textit{power sum symmetric functions} of degree $n$, $$\{p_{\lambda}=p_{\lambda_1}\cdots p_{\lambda_k}: \lambda=(\lambda_1,\ldots,\lambda_k)\vdash  n\}.$$

A concrete realization of $\Lambda$ is obtained by embedding $\Lambda=\mathbb{Q}[p_1,p_2,\ldots]$ in $\mathbb{Q}[[x_1,x_2,\ldots]]$, i.e. the ring of formal power series in countably many commuting indeterminates $\{x_1,x_2,\ldots\}$, under the identification (extended multiplicatively)
$$p_i \longmapsto \sum_{j\geq 1}x_j^i.$$
Then we can think of symmetric functions as being formal power series $f$ in the $x$ variables with the property that $f(x_{\pi(1)},x_{\pi(2)},\ldots)=f(x_{1},x_{2},\ldots)$ for every permutation $\pi$ of the positive integers $\mathbb{N}$. It is with this viewpoint that we will define a very important class of symmetric functions next. 

\subsection{Schur functions}
We start by defining the skew Schur functions combinatorially. 
\begin{Definition} Given a skew shape $\lambda / \mu$, the \emph{skew Schur function} of \emph{shape} $\lambda / \mu$, $s_{\lambda / \mu}$, is the formal power series
\begin{eqnarray*}
s_{\lambda / \mu}=\displaystyle\sum_{T\in SSYT(\lambda/\mu)}x^{T}.
\end{eqnarray*}
If $\mu=\varnothing$, then $\lambda/\mu=\lambda$, and we call $s_{\lambda}$ the \emph{Schur function} of shape $\lambda$.
\end{Definition}
Though not evident from the definition, skew Schur functions are actually symmetric functions and the elements of the set $\{s_{\lambda}:\lambda\vdash n\}$ form a basis for $\Lambda^n$.
Now we can equip this space with an inner product $\langle , \rangle_{\Lambda^{n}}$, called the \textit{Hall inner product}. It is defined by setting $\langle s_{\lambda},s_{\mu}\rangle_{\Lambda^{n}} = \delta_{\lambda\mu}$, where $\delta_{\lambda\mu}=1$ if $\lambda=\mu$ and $0$ otherwise, and then defining the inner product for any $f,g \in \Lambda^{n}$ by linear extension. One can extend this to an inner product on $\Lambda$, in which case we will refer to it as $\langle ,\rangle_{\Lambda}$. We will need the following very fundamental property of skew Schur functions, given partitions $\lambda, \mu$ and $\nu$. 
\begin{eqnarray*}
\langle s_{\mu}s_{\nu},s_{\lambda}\rangle_{\Lambda}=\langle s_{\nu},s_{\lambda/\mu}\rangle_{\Lambda}
\end{eqnarray*}

There exists a combinatorial rule to multiply two Schur functions and express the result in the Schur basis by counting SSYTs satisfying certain constraints, commonly called the \textit{Littlewood-Richardson rule}. Thus, given partitions $\mu$ and $\nu$, we have an expansion as follows
\begin{eqnarray*}
s_{\mu}s_{\nu}=\displaystyle \sum_{\lambda}c_{\mu ,\nu}^{\lambda}s_{\lambda}
\end{eqnarray*}
where the sum is over all $\lambda$ such that $\mu$ is contained in $\lambda$. Here the $c_{\mu,\nu}^{\lambda}$ turn out to be non-negative integers that can be computed combinatorially, and are called the \textit{Littlewood-Richardson coefficients}. In terms of the inner product on $\Lambda$ the above is equivalent to $\langle s_{\mu}s_{\nu},s_{\lambda}\rangle_{\Lambda}=\langle s_{\nu},s_{\lambda/\mu}\rangle_{\Lambda}=c_{\mu,\nu}^{\lambda}$. 

We will only require special cases of the Littlewood-Richardson rule which describe the multiplication of a Schur function with a Schur function indexed by shape of one row or one column. These cases are collectively called the \textit{Pieri rule} but before we state the rule we need to describe certain skew shapes. A skew shape $\lambda / \mu$ is called a \textit{horizontal strip} if it does not contain boxes in the same column, and is called a \textit{vertical strip} if it does not contain boxes in the same row.
\begin{Theorem}[Pieri rule]\label{theorem: Pieri rules}
If $\mu$ is a partition, then
\begin{eqnarray*}
s_{\mu}s_{(n)}&=&\displaystyle\sum_{\substack{\nu \vdash |\mu |+n\\ \nu / \mu = \text{horizontal strip of size n}}}s_{\nu}\\
s_{\mu}s_{(1^n)}&=&\displaystyle\sum_{\substack{\nu \vdash |\mu |+n\\ \nu / \mu = \text{vertical strip of size n}}}s_{\nu}.
\end{eqnarray*}
\end{Theorem}
\subsection{The Kronecker product of Schur functions}
In this section we will outline how the Kronecker coefficients arise in the representation theory of the symmetric group. Given $\mu\vdash n$, let $V^{\mu}$ denote the irreducible representation of $\mathfrak{S}_n$ indexed by $\mu$, whose dimension equals $f_{\mu}$, and let the corresponding character be denoted by $\chi_{\mu}$. Then the pointwise product $\chi_{\mu}\chi_{\nu}$ is the character of the $\mathfrak{S}_n$-representation $V^{\mu}\otimes V^{\nu}$. Let $g_{\mu\nu}^{\lambda}$ be the multiplicity of $V^{\lambda}$ in $V^{\mu}\otimes V^{\nu}$. That is, $g_{\mu\nu}^{\lambda}=\langle \chi_{\mu}\chi_{\nu} , \chi_{\lambda}\rangle_{CF^n} $ where $\langle ,\rangle _{CF^n}$ denotes the standard inner product on the space of class functions $CF^n$ of $\mathfrak{S}_n$. We will now explain how this interpretation ties in with the Kronecker product for Schur functions that we will define soon.

The \textit{Kronecker product}, $*$, on $\Lambda$ is defined implicitly by defining it on the basis of power sum symmetric functions as follows
\begin{eqnarray*}
\frac{p_{\lambda}}{z_{\lambda}}*\frac{p_{\mu}}{z_{\mu}}=\delta_{\lambda\mu}\frac{p_{\lambda}}{z_{\lambda}},
\end{eqnarray*}
and then extending it linearly. Here $z_{\lambda}$ denotes the number of permutations in $\mathfrak{S}_{|\lambda|}$ commuting with a fixed permutation of cycle type $\lambda$ where the \textit{cycle type} of a permutation $\sigma$ is the partition obtained by ordering the cycle lengths occurring in the cycle decomposition of $\sigma$ in weakly decreasing order.

With this definition, it turns out that 
\begin{eqnarray}\label{eqn: Kronecker both interpretations}
g_{\mu\nu}^{\lambda}=\langle\chi_{\mu}\chi_{\nu}, \chi_{\lambda}\rangle_{CF^n}=\langle s_{\mu}*s_{\nu}, s_{\lambda}\rangle_{\Lambda^n}
\end{eqnarray}
where $\lambda$, $\mu$ and $\nu$ are partitions of the same size $n$. The Kronecker product also satisfies the following
\begin{eqnarray*}
s_{\mu}*s_{\nu}=s_{\nu}*s_{\mu} \text{ and } s_{\mu}*s_{\nu}=s_{\nu^t}*s_{\mu^t}.
\end{eqnarray*}
Moreover, if $\mu,\nu \vdash n$ then
$g_{\mu\nu}^{(n)}=g_{\mu\nu^t}^{(1^{n})}=\delta_{\mu\nu}$.
\begin{Remark}
Since we will never be using the inner product on the space of class functions, $\langle ,\rangle$ will always mean the Hall inner product $\langle ,\rangle_{\Lambda}$ from now on.
\end{Remark}
Before we recall the relevant results on Kronecker products, we will establish some notation that we will stick to throughout. Given a positive integer $n$, let 
\begin{eqnarray*}
P_{n}&=&\{\lambda \vdash 2n: l(\lambda)\leq 4\text{ and } \lambda \text{ has either all parts even or $l(\lambda)=4$ and all parts odd})\},\\
Q_{n}&=&\{\lambda \vdash 2n:  l(\lambda) \leq 4 \text{ and exactly two of } \lambda_{i} \text{ are odd})\}.
\end{eqnarray*}
This given, let 
\begin{eqnarray*}
P=\bigcup_{n\geq 0}P_{n}\hspace{2mm}\text{ and }Q=\bigcup_{n\geq 0}Q_{n}\hspace{2mm},
\end{eqnarray*}
and it is amply clear that $P\cup Q$ is the set of all partitions of even size and length at most $4$.

We will also be needing the Knuth bracket for giving truth values to statements. 
\begin{eqnarray*}
((S))=\left\lbrace\begin{array}{ll} 1 &  S \text{ is a true statement}\\ 0  & \text{otherwise}\end{array} \right.
\end{eqnarray*}
Now we are in a position to state the results of interest to us. The computation of $s_{(n,n)}*s_{(n,n)}$ is one such result. This computation originally arose out of solving a mathematical physics problem related to resolving the interference of 4 qubits \cite{W}. It appeared first in \cite{GWXZ} in the form as shown below. It was proven  again in \cite{BvWZ}. The result states the following.
\begin{Theorem}\cite[Theorem I.6]{GWXZ}\label{nnprod}
Given a positive integer $n$, 
\begin{eqnarray}
s_{(n,n)}*s_{(n,n)}=\sum_{\lambda \in P_n}s_{\lambda}.
\end{eqnarray}
\end{Theorem}
This characterization is different from earlier characterizations as it explicitly states which partitions have non-zero coefficients and further establishes that the coefficients are all either 0 or 1 without giving a combinatorial rule. 

Using the result of \cite{GWXZ} as inspiration, a characterization of the Kronecker product of $s_{(n,n)}*s_{(n+k,n-k)}$ for $k \geq 0$  was obtained in \cite{BvWZ}. Since we do not need the full strength of their result, we will just state the $k=1$ case.
\begin{Theorem}\cite[Corollary 3.6]{BvWZ}\label{nplusminusprod}
Given a positive integer $n$, 
\begin{eqnarray}
s_{(n+1,n-1)}*s_{(n,n)}=\sum_{\lambda \in Q_n}s_{\lambda}.
\end{eqnarray}
\end{Theorem}

A result of Littlewood that we will frequently use, and which simplifies calculations at many places is the following \cite{L}.
\begin{Theorem}[Littlewood]\label{theorem: Littlewood}
Let $\alpha,\beta$ and $\gamma$ be partitions such that $|\alpha|+|\beta|=|\gamma|$. Then,
\begin{equation*}
(s_{\alpha}s_{\beta})*s_{\gamma}=\displaystyle \sum_{\delta \vdash |\beta|}\sum_{\eta \vdash |\alpha|} c_{\eta,\delta}^{\gamma} (s_{\eta}*s_{\alpha})(s_{\delta}*s_{\beta}) 
\end{equation*} 
where $c_{\eta,\delta}^{\gamma}$ are the Littlewood-Richardson coefficients.
\end{Theorem}

Using this identity of Littlewood in conjunction with Theorem \ref{nnprod}, one can prove the following corollary, present in the following form in \cite{BvWZ}.
\begin{Corollary}\cite[Corollary 4.1]{BvWZ}\label{nminusprod}
Given a positive integer $n$,
\begin{eqnarray*}
s_{(n,n-1)}*s_{(n,n-1)}=\sum_{\substack{\lambda \vdash 2n-1\\l(\lambda)\leq 4}}s_{\lambda}.
\end{eqnarray*}
\end{Corollary}
We will need one final result which, given partitions $\mu$ and $\nu$, helps in identifying certain partitions $\lambda$ for which $g_{\mu\nu}^{\lambda}=0$. Below, $\mu \cap \nu$ denotes the partition obtained by intersecting the corresponding Ferrers diagrams once their top left corners are aligned. Clausen and Meier \cite{CM} and Dvir \cite{D} proved the following theorem.
\begin{Theorem}\label{theorem: dvir, clausen and meier}
Let $\mu$, $\nu$ be partitions of $n$. Then 
\begin{eqnarray*}
\text{max }\{\lambda_{1}:g_{\mu\nu}^{\lambda} \neq 0 \text{ for some } \lambda =(\lambda_{1},\ldots , \lambda_{l(\lambda)})\}=|\mu \cap \nu |,\\
\text{max }\{l(\lambda):g_{\mu\nu}^{\lambda} \neq 0\text{ for some } \lambda =(\lambda_{1},\ldots ,\lambda_{l(\lambda)})\}=|\mu \cap \nu^t|.
\end{eqnarray*}
\end{Theorem}
The import of this theorem can be gauged by the fact that it already implies that if $\mu$ and $\nu$ are partitions each with at most two rows, then $g_{\mu\nu}^{\lambda}=0$ for all $\lambda$ such that $l(\lambda) \geq 5$.

\section{The Kronecker coefficient $g_{(n,n-1,1)(n,n)}^{\theta}$ in $s_{(n,n-1,1)}*s_{(n,n)}$ where $n \geq 2$}\label{section:first case}
We will now give an explicit characterization of the Kronecker product of 
$s_{(n,n-1,1)}$ and $s_{(n,n)}$. Observe that the Pieri rule (Theorem \ref{theorem: Pieri rules}) implies 
\begin{equation} s_{(n,n-1,1)} = s_{(n,n-1)}s_{(1)}-s_{(n,n)}-s_{(n+1,n-1)} .\label{firststep}\end{equation}
Since we are interested in computing the coefficients $g_{(n,n-1,1)(n,n)}^{\theta}$ where $\theta \vdash 2n$, we compute $\langle s_{(n,n-1,1)}*s_{(n,n)},s_{\theta} \rangle$ by \eqref{eqn: Kronecker both interpretations}. Using $\eqref{firststep}$, we obtain
 \begin{eqnarray} \langle s_{(n,n-1,1)}*s_{(n,n)},s_{\theta} \rangle &=& \langle  (s_{(n,n-1)}s_{(1)})*s_{(n,n)},s_{\theta} \rangle \nonumber \\ & &- \langle (s_{(n,n)}+s_{(n+1,n-1)})*s_{(n,n)},s_{\theta} \rangle .\label{secondstep}\end{eqnarray}
 We will evaluate the inner products appearing on the right hand side of $\eqref{secondstep}$ individually. The use of Theorem $\ref{theorem: Littlewood}$ implies
\begin{eqnarray}
(s_{(n,n-1)}s_{(1)})*s_{(n,n)} &=& \displaystyle \sum_{\delta \vdash 1}\sum_{\eta \vdash 2n-1} c_{\eta,\delta}^{(n,n)} (s_{\eta}*s_{(n,n-1)})(s_{\delta}*s_{(1)}) \nonumber \\ &=& \displaystyle \sum_{\eta \vdash 2n-1} c_{\eta,(1)}^{(n,n)} (s_{\eta}*s_{(n,n-1)})(s_{(1)}*s_{(1)}).
\end{eqnarray}
The Pieri rule yields that $c_{\eta,(1)}^{(n,n)} \neq 0$ if and only if $\eta = (n,n-1)$, in which case $c_{(n,n-1),(1)}^{(n,n)} = 1$. Since $s_{(1)}*s_{(1)} = s_{(1)}$, we conclude that 
\begin{eqnarray}
(s_{(n,n-1)}s_{(1)})*s_{(n,n)}&=&s_{(1)}(s_{(n,n-1)}*s_{(n,n-1)}).
\end{eqnarray}
This reduces $\eqref{secondstep}$ to 
\begin{eqnarray} \langle s_{(n,n-1,1)}*s_{(n,n)},s_{\theta} \rangle &=& \langle  s_{(1)}(s_{(n,n-1)}*s_{(n,n-1)}),s_{\theta} \rangle \nonumber \\ & &- \langle (s_{(n,n)}+s_{(n+1,n-1)})*s_{(n,n)},s_{\theta} \rangle \nonumber \\ 
 &=&\langle  s_{(n,n-1)}*s_{(n,n-1)},s_{\theta /(1)} \rangle \nonumber \\ & &- \langle (s_{(n,n)}+s_{(n+1,n-1)})*s_{(n,n)},s_{\theta} \rangle \nonumber \\ &=& \langle \displaystyle\sum_{\substack{\lambda \vdash 2n-1\\ l(\lambda) \leq 4}}s_{\lambda}, s_{\theta/(1)} \rangle -\langle \displaystyle\sum_{\substack{\lambda \vdash 2n \\ l(\lambda) \leq 4}}s_{\lambda},s_{\theta} \rangle . \label{firstsequence}\end{eqnarray}

In arriving at the last step in the above sequence, we have made use of Corollary \ref{nminusprod}, Theorem \ref{nnprod} and Theorem \ref{nplusminusprod}. Notice $\langle \displaystyle\sum_{\lambda \vdash 2n, l(\lambda) \leq 4}s_{\lambda},s_{\theta} \rangle$ is $1$ if $l(\theta) \leq 4$ and $0$ otherwise. So we will focus on evaluating $\langle \displaystyle\sum_{\lambda \vdash 2n-1,l(\lambda) \leq 4}s_{\lambda}, s_{\theta/(1)} \rangle$. The Pieri rule implies 
\begin{eqnarray}
s_{\theta/(1)}&=&\displaystyle\sum_{\theta^{-} \prec \theta}s_{\theta^{-}}.
\label{goingdown}
\end{eqnarray}
Note the crucial fact that the number of terms appearing on the right hand side of \eqref{goingdown} is equal to the number of distinct parts in the partition $\theta$, i.e. $d_{\theta}$.

If $\langle s_{(n,n-1,1)}*s_{(n,n)},s_{\theta} \rangle \neq 0$, then we must have that $l(\theta) \leq 5$ by Theorem \ref{theorem: dvir, clausen and meier}, as $|(n,n-1,1) \cap (n,n)^{t}| \leq 5$. We will carry out the rest of the computation in cases depending on the length of the partition $\theta$.

\subsection{Case I: $\mathbf{l(\theta)=5}$}
If $l(\theta)=5$, but $\theta_{5} \geq 2$, then $s_{\theta/(1)}$ is sum of terms of the form $s_{\gamma}$ with $l(\gamma)=5$. The right hand side of $\eqref{firstsequence}$ clearly implies that the coefficient of $s_{\theta}$ in $s_{(n,n-1,1)}*s_{(n,n)}$ is $0$ in this instance. If $\theta_{5}=1$, then $s_{\theta/(1)} = s_{\theta^{'}} + $ sum of terms of the form $s_{\gamma}$ where $l(\gamma)=5$. This in turn means that $\langle \displaystyle\sum_{\substack{\lambda \vdash 2n-1\\ l(\lambda) \leq 4}}s_{\lambda}, s_{\theta/(1)} \rangle = 1$. Thus, if $l(\theta)=5$,
\begin{eqnarray*}
\langle s_{(n,n-1,1)}*s_{(n,n)},s_{\theta} \rangle =  \left\lbrace \begin{array}{ll}1 & \theta_{5}=1\\0 & \text{otherwise.}\end{array}\right.
\end{eqnarray*}

\subsection{Case II: $\mathbf{l(\theta)\leq 4}$}
We know that if $l(\theta) \leq 4$, then 
$\langle \displaystyle\sum_{\lambda\vdash 2n, l(\lambda)\leq 4}s_{\lambda},s_{\theta} \rangle = 1$. The following computation helps us in finishing this case.
\begin{eqnarray}
\langle \displaystyle\sum_{\substack{\lambda \vdash 2n-1\\ l(\lambda) \leq 4}}s_{\lambda}, s_{\theta/(1)} \rangle &=& \langle \displaystyle\sum_{\substack{\lambda \vdash 2n-1\\ l(\lambda) \leq 4}}s_{\lambda},\displaystyle\sum_{\theta^{-} \prec \theta}s_{\theta^{-}}\rangle \nonumber \\ &=& d_{\theta}.
\end{eqnarray} 
Thus, using \eqref{firstsequence}, we get that  for $l(\theta)\leq 4$,
\begin{eqnarray*}
\langle s_{(n,n-1,1)}*s_{(n,n)},s_{\theta} \rangle =  d_{\theta} - 1 .
\end{eqnarray*}

On collecting the results of the two cases together, we obtain the following description.
\begin{mdframed}
\begin{eqnarray*}
g_{(n,n-1,1)(n,n)}^{\theta}=  \left\lbrace \begin{array}{ll}1 & l(\theta)=5,\text{ }\theta_{5}=1\\d_{\theta}-1 & l(\theta)\leq 4\\0 &\text{otherwise}\end{array}\right.
\end{eqnarray*}
\end{mdframed}

\begin{Example}
Consider the computation of $s_{(4,3,1)}*s_{(4,4)}$. Then
\begin{eqnarray*}
s_{(4,3,1)}*s_{(4,4)}=&&s_{(2,2,2,1,1)}+s_{(3,2,1,1,1)}+2s_{(3,2,2,1)}+s_{(3,3,1,1)}+s_{(3,3,2)}\\&&+s_{(4,1,1,1,1)}+2s_{(4,2,1,1)}+s_{(4,2,2)}+2s_{(4,3,1)}+s_{(5,1,1,1)}\\&&+2s_{(5,2,1)}+s_{(5,3)}+s_{(6,1,1)}+s_{(6,2)}+s_{(7,1)}.
\end{eqnarray*}
\end{Example}

\section{The Kronecker coefficient $g_{(n-1,n-1,1)(n,n-1)}^{\theta}$ in $s_{(n-1,n-1,1)}*s_{(n,n-1)}$ where $n\geq 2$}\label{section: second case}
In the same vein as the previous case, we can explicitly compute the Kronecker product of $s_{(n-1,n-1,1)}$ and $s_{(n,n-1)}$. Again, the Pieri rule implies that  
\begin{eqnarray}
s_{(n-1,n-1,1)} &=& s_{(1)}s_{(n-1,n-1)}-s_{(n,n-1)}.
\label{thirdstep}
\end{eqnarray}

An application of Theorem $\ref{theorem: Littlewood}$ gives 
\begin{eqnarray}
(s_{(1)}s_{(n-1,n-1)})*s_{(n,n-1)}&=&\displaystyle \sum_{\delta \vdash 1}\sum_{\eta \vdash 2n-2} c_{\eta,\delta}^{(n,n-1)} (s_{\eta}*s_{(n-1,n-1)})(s_{\delta}*s_{(1)}) \nonumber \\ &=&\displaystyle\sum_{\eta \vdash 2n-2} c_{\eta,(1)}^{(n,n-1)} (s_{\eta}*s_{(n-1,n-1)})(s_{(1)}*s_{(1)}).
\end{eqnarray}

The Pieri rule dictates that the only cases where $c_{\eta,(1)}^{(n,n-1)} \neq 0$ are when $\eta = (n,n-2)$ or $\eta = (n-1,n-1)$ and in both cases $c_{\eta,(1)}^{(n,n-1)} = 1$. Thus
\begin{eqnarray}
(s_{(1)}s_{(n-1,n-1)})*s_{(n,n-1)}&=& s_{(1)}(s_{(n,n-2)}*s_{(n-1,n-1)})\nonumber \\& & +\textbf{ }s_{(1)}(s_{(n-1,n-1)}*s_{(n-1,n-1)}).
\label{fourthstep}
\end{eqnarray}

If $\theta \vdash 2n-1$, $\eqref{thirdstep}$ and $\eqref{fourthstep}$ together bring us to 
\begin{eqnarray}
\langle s_{(n-1,n-1,1)}*s_{(n,n-1)},s_{\theta} \rangle &=& \langle s_{(1)}((s_{(n,n-2)}+s_{(n-1,n-1)})*s_{(n-1,n-1)}),s_{\theta}\rangle \nonumber \\& &-\textbf{ }\langle s_{(n,n-1)}*s_{(n,n-1)},s_{\theta} \rangle \nonumber \\ &=& \langle (s_{(n,n-2)}+s_{(n-1,n-1)})*s_{(n-1,n-1)},s_{\theta / (1)}\rangle \nonumber \\& &-\textbf{ }\langle s_{(n,n-1)}*s_{(n,n-1)},s_{\theta} \rangle \nonumber \\ &=&\langle \displaystyle\sum_{\substack {\lambda \vdash 2n-2\\l(\lambda) \leq 4}}s_{\lambda},s_{\theta/(1)}\rangle -\langle \displaystyle\sum_{\substack{\lambda \vdash 2n-1\\l(\lambda) \leq 4}}s_{\lambda},s_{\theta}\rangle . 
\label{secondcaseequation}
\end{eqnarray}

Now one can easily check that this gives the same characterization as the one obtained from $s_{(n,n-1,1)}*s_{(n,n)}$ and the argument is essentially the same, except that we use \eqref{secondcaseequation} instead of \eqref{firstsequence}. Hence we obtain the following.
\begin{mdframed}
\begin{eqnarray*}
 g_{(n-1,n-1,1)(n,n-1)}^{\theta} = \left \{\begin{array}{ll}1 & l(\theta)=5,\textbf{ } \theta_{5}=1\\ d_{\theta}-1 & l(\theta) \leq 4 \\0 & \text{otherwise}\end{array}\right. 
\end{eqnarray*}
\end{mdframed}

\begin{Example}
We will compute $s_{(3,3,1)}*s_{(4,3)}$.
\begin{eqnarray*}
s_{(3,3,1)}*s_{(4,3)}=&&s_{(2,2,1,1,1)}+s_{(2,2,2,1)}+s_{(3,1,1,1,1)}+2s_{(3,2,1,1)}+s_{(3,2,2)}+s_{(3,3,1)}\\&&+s_{(4,1,1,1)}+2s_{(4,2,1)}+s_{(4,3)}+s_{(5,1,1)}+s_{(5,2)}+s_{(6,1)}.
\end{eqnarray*}
\end{Example}

\section{The Kronecker coefficient $g_{(n-1,n-1,2)(n,n)}^{\theta}$ in $s_{(n-1,n-1,2)}*s_{(n,n)}$ where $n \geq 3$}\label{section: third case}
Before we give the Kronecker coefficients occurring in the product $s_{(n-1,n-1,2)}*s_{(n,n)}$, we will make one remark about our notation. Henceforth, the statement `$\lambda \in P$' is considered to be equivalent to `$\lambda' \in P$', and an analogous statement holds for a statement like `$\lambda \in Q$'. For example, consider $\lambda=(5,3,3,1,1)$. Then, even though $\lambda$ has 5 parts, we say $((\lambda \in P))$ evaluates to $1$ because $\lambda'=(5,3,3,1)$ has all 4 parts odd, and thus belongs to $P$.

Next we give a formula for the Kronecker coefficient $g_{(n-1,n-1,2)(n,n)}^{\theta}$ where $\theta\vdash 2n$ and follow it up with an example. Its proof can be found in \cite[Section 2.3]{Tewari}.
\begin{mdframed}
\begin{eqnarray*}g_{(n-1,n-1,2)(n,n)}^{\theta}  = \left \{\begin{array}{ll}((\theta \in P)) & l(\theta)=6\text{ and } \theta_{5}=\theta_{6}=1\\ ((\theta \in Q)) & l(\theta)=5,\text{ }\theta_{5}=2\\O_{\theta}^{'}-((\theta_{4}=1)) & l(\theta)=5, E_{\theta} =1\text{ and }\theta_{5}=1\\ E_{\theta}^{'} & l(\theta)=5, O_{\theta} = 1\text{ and }\theta_{5}=1\\ 1-d_{\theta}+\displaystyle\binom{d_{\theta}}{2} & l(\theta)\leq 4,\text{ }\theta \in P\\ 1-d_{\theta}+d_{\theta,2}+O_{\theta}^{'}E_{\theta}^{'}+((E_{\theta}=2)) & l(\theta)\leq 4,\text{ }\theta \in Q\\0 & \text{otherwise}\end{array}\right. \end{eqnarray*} \end{mdframed}

\begin{Example}
Consider the product $s_{(7,7,2)}*s_{(8,8)}$, and three partitions $\alpha=(5,5,3,1,1,1)$, $\beta=(6,4,3,2,1)$ and $\gamma= (7,5,2,2)$. 

Note that $l(\alpha)=6$ and $\alpha_{5}=\alpha_{6}=1$ as well. Since $\alpha'=(5,5,3,1) \in P$, we obtain 
$g_{(7,7,2)(8,8)}^{(5,5,3,1,1,1)}=1$.

Consider the case of $\beta$ now. We have $l(\beta)=5$ and $\beta_{5}=1$. Since $\beta'$ has exactly 1 odd part, we have $O_{\beta}=1$. Thus, the above characterization allows us to obtain $g_{(7,7,2)(8,8)}^{(6,4,3,2,1)}=E_{(6,4,3,2,1)}^{'}$. Since $\beta'$ has exactly 3 distinct even parts, we have $g_{(7,7,2)(8,8)}^{(6,4,3,2,1)}=3$.

Turning our attention to $\gamma$, we see that $l(\gamma)=4$ and $\gamma \in Q$. We have $d_{\gamma}=3$, $d_{\gamma,2}=3$,  $O_{\gamma}^{'}=2$, $E_{\gamma}^{'}=1$ and $E_{\gamma}=2$. Thus, we get $g_{(7,7,2)(8,8)}^{(7,5,2,2)}=4$.
\end{Example}

\section{The Kronecker coefficient $g_{(n-1,n-1,1,1)(n,n)}^{\theta}$ in $s_{(n-1,n-1,1,1)}*s_{(n,n)}$ where $n\geq 2$}\label{section: fourth case}
First we give a formula for the Kronecker coefficient $g_{(n-1,n-1,1,1)(n,n)}^{\theta}$ where $\theta\vdash 2n$ and then an example. Its proof can be found in \cite[Section 2.4]{Tewari}.
\begin{mdframed}
\begin{eqnarray*}g_{(n-1,n-1,1,1)(n,n)}^{\theta} &=&  \left \{\begin{array}{ll}((\theta\in Q)) & l(\theta)=6\text{ and }\theta_{5}=\theta_{6}=1\\1 & l(\theta)=5, \theta_{5}=2\text{ and }\theta \in P\\O_{\theta}^{'} & l(\theta)=5, E_{\theta}=1\text{ and }\theta_{5}=1\\E_{\theta}^{'}-((\theta_{4}=1)) & l(\theta)=5, O_{\theta}=1 \text{ and  }\theta_{5}=1\\d_{\theta,2}-d_{\theta} +\displaystyle\binom{d_{\theta}}{2}+R_{\theta} & l(\theta)\leq 4,\text{ }\theta \in P\\1-d_{\theta}+((E_{\theta}^{'}=2))\\+((O_{\theta}^{'}=2))+O_{\theta}^{'}E_{\theta}^{'} & l(\theta)\leq 4,\text{ }\theta \in Q\\0 & \text{otherwise} \end{array} \right. \end{eqnarray*}\end{mdframed}
\begin{Example}
Consider the Kronecker product $s_{(7,7,1,1)}*s_{(8,8)}$, and three partitions $\alpha=(5,5,3,1,1,1)$, $\beta=(6,4,3,2,1)$ and $\gamma= (7,5,2,2)$.

Notice that even though $\alpha_{5}=\alpha_{6}=1$, the partition $\alpha'=(5,5,3,1)\notin Q$. Thus $g_{(7,7,1,1)(8,8)}^{(5,5,3,1,1,1)}=0$.

Consider the case of $\beta$ now. We have $l(\beta)=5$ and $\beta_{5}=1$. Since $\beta'=(6,4,3,2)$, we have $O_{\beta}=1$. The characterization above gives us $g_{(7,7,1,1)(8,8)}^{(6,4,3,2,1)}= E_{(6,4,3,2,1)}^{'}-((\beta_{4}=1))$. Since $\beta'$ has 3 distinct even parts and $\beta_{4} \neq 1$, we obtain $g_{(7,7,1,1)(8,8)}^{(6,4,3,2,1)}=3$.

As far as $\gamma$ is concerned, we have $l(\gamma)=4$ and $\gamma \in Q$. We have $d_{\gamma}=3$, $O_{\gamma}^{'}=2$ and $E_{\gamma}^{'}=1$. Thus $g_{(7,7,1,1)(8,8)}^{(7,5,2,2)}= 1-3+0+1+2=1$.
\end{Example}

\section{The Kronecker coefficient $g_{(n,n,1)(n,n,1)}^{\theta}$ in $s_{(n,n,1)}*s_{(n,n,1)}$ where $n\geq 2$}\label{section: fifth case}
In this section we will compute the Kronecker product $s_{(n,n,1)}*s_{(n,n,1)}$. Before we begin our calculations, we need to introduce certain statistics on partitions, the purpose of which is spelt out in detail in \cite[Section 2.5]{Tewari}. In essence, we need these statistics to deduce the relation between the number of distinct parts in a partition $\theta$, and that in a partition $\theta^{-}\prec\theta$. 

Fix an alphabet $X=\{0,1,2\}$. We will associate a string $\sigma$ of length $l(\theta)+1$ to a partition $\theta$. For $1\leq i \leq l(\theta)$, define
\begin{eqnarray*}
\sigma_{i}= \left\lbrace \begin{array}{ll}0 & (\theta_{i}=\theta_{i+1})\\1 & (\theta_{i}-\theta_{i+1}=1)\\2 & (\theta_{i}-\theta_{i+1} \geq 2).\end{array}\right.
\end{eqnarray*}
Here we are assuming that when $i=l(\theta)$, then $\theta_{i+1}=0$. For the sake of convenience, define $\sigma_{0}=\sigma_{1}$. Once $\sigma$ has been found, define the following sets.
\begin{eqnarray*}
A_{\theta,1}&=&\{i:1\leq i\leq l(\theta), \text{ }\sigma_{i}=1 \text{ and }\sigma_{i-1}=0\}\nonumber\\
A_{\theta,2}&=&\{i:1\leq i\leq l(\theta), \text{ }\sigma_{i}=2 \text{ and }\sigma_{i-1}=0\}\nonumber\\
B_{\theta,1}&=&\{i:1\leq i\leq l(\theta), \text{ }\sigma_{i}=1 \text{ and }\sigma_{i-1}\neq 0 \}\nonumber\\
B_{\theta,2}&=&\{i:1\leq i\leq l(\theta), \text{ }\sigma_{i}=2 \text{ and }\sigma_{i-1}\neq 0\}
\end{eqnarray*}
We will use $a_{\theta,1}$, $a_{\theta,2}$, $b_{\theta,1}$ and $b_{\theta,2}$ to denote the cardinalities of the sets $A_{\theta,1}$, $A_{\theta,2}$, $B_{\theta,1}$ and $B_{\theta,2}$ respectively.
Now we can give a description for the Kronecker coefficients occurring in $s_{(n,n,1)}*s_{(n,n,1)}$.
\begin{mdframed}\begin{eqnarray*}
g_{(n,n,1)(n,n,1)}^{\theta} = \left\lbrace \begin{array}{ll}1 & l(\theta)=6,\theta_{6}=\theta_{5}=1\\2d_{\theta}-3+((\theta\in P)) & l(\theta)=5,\theta_{5}=\theta_{4}=1\\2d_{\theta}-4+((\theta\in P)) & l(\theta)=5,\theta_{4} \geq 2,\theta_{5}=1\\1 & l(\theta)=5,\theta_{5}=2\\(d_{\theta}-1)^2+1-b_{\theta ,1}+a_{\theta ,2} & l(\theta)=4\\(d_{\theta}-1)^2+1-b_{\theta ,1}+a_{\theta ,2} & l(\theta)=3 \text{ and }\theta \\&\text{has exactly 1 odd part} \\(d_{\theta}-1)^2-b_{\theta ,1}+a_{\theta ,2} & l(\theta)=3 \text{ and }\theta \\&\text{has all parts odd}\\2-b_{\theta ,1}+a_{\theta ,2} & l(\theta)=2\\1& l(\theta)=1\\0 & \text{otherwise}\end{array}\right.
\end{eqnarray*}
\end{mdframed}
Let us now consider an example.
\begin{Example}
Consider the computation of $s_{(8,8,1)}*s_{(8,8,1)}$. Let $\alpha=(6,5,3,2,1)$, $\beta= (8,6,2,1)$ and $\gamma = (7,5,5)$ be three partitions.

Consider first the case of $\alpha$. We can see that $l(\alpha)=d_{\alpha}=5$, $\alpha_{5}=1$ and $\alpha_{4} \geq 2$. Note also that $\alpha'=(6,5,3,2)$ does not belong to $P$. Thus $ g_{(8,8,1)(8,8,1)}^{(6,5,3,2,1)}=2\times 5-4=6$.

Next, consider $\beta$. We have $l(\beta)=d_{\beta}=4$. The string $\sigma$ associated with $(8,6,2,1)$ is $22211$. This immediately yields $a_{\beta,2}=0$ and $b_{\beta,1}=2$. This implies $g_{(8,8,1)(8,8,1)}^{(8,6,2,1)}=(4-1)^2+1-2+0=8$.

Finally, consider $\gamma$. We can see that $l(\gamma)=3$, $d_{\gamma}=2$ and the string $\sigma$ associated with $\gamma=(7,5,5)$ is $2202$. Thus, we have $a_{\gamma,2}=1$ and $b_{\gamma,1}=0$. This gives $g_{(8,8,1)(8,8,1)}^{(7,5,5)} = (2-1)^2-0+1=2$.
\end{Example}

\section{Combinatorial implications}\label{section: enumerative applications}
A natural question one can ask is how many SYTs are there of fixed size $n$ if one imposes the constraint that the number of parts of $\lambda \vdash n$ is bounded above by some fixed positive integer $k$. This means we are interested in the sum
\begin{eqnarray*}
\tau_{k}(n)=\displaystyle\sum_{\substack{\lambda\vdash n\\l(\lambda)\leq k}}f_{\lambda}.
\end{eqnarray*}
This is also a well studied question as is evident from \cite{BFK,BG,G,R}. The expressions for $\tau_{k}(n)$ are unwieldy when $k$ is large. But for relatively small values of $k$, these expressions are more succinct than what one would expect from the hooklength formula. Regev \cite{R} found the following closed form expressions for $\tau_{2}(n)$ and $\tau_{3}(n)$.
\begin{eqnarray}\label{eqn:motzkin}
\tau_{2}(n)=\displaystyle\binom{n}{\lfloor\frac{n}{2}\rfloor},\hspace{3mm}\tau_{3}(n)=\sum_{i\geq 0}\frac{1}{i+1}\binom{n}{2i}\binom{2i}{i}.
\end{eqnarray}
Note that $\tau_{3}(n)$ is actually the Motzkin number $M_{n}$. Gessel \cite{G} found an expression for $\tau_{4}(n)$ while Gouyou-Beauchamps \cite{Gb} found an expression for both $\tau_{4}(n)$ and $\tau_{5}(n)$. The expressions were
\begin{eqnarray}\label{eqn:catalan}
\tau_{4}(n)=C_{\lfloor\frac{n+1}{2}\rfloor}C_{\lceil\frac{n+1}{2}\rceil}, \hspace{3mm}\tau_{5}(n)=6\displaystyle\sum_{i=0}^{\lfloor\frac{n}{2}\rfloor}\binom{n}{2i}C_{i}\frac{(2i+2)!}{(i+2)!(i+3)!},
\end{eqnarray}
where $C_{i}=\frac{1}{i+1}\binom{2i}{i}$ is the $i$-th Catalan number. For a comprehensive list of the many interpretations of both Motzkin and Catalan numbers, the reader should refer to \cite{St}.

In this section, we will use our results on Kronecker coefficients to prove a result similar in nature to the ones mentioned above. Given a positive integer $n$, consider the set $L_n$ defined as follows.
\begin{eqnarray*}
L_n=\{\lambda\vdash n: l(\lambda)=5, \lambda_5 =1\}
\end{eqnarray*}
In Theorem \ref{theorem:main combinatorial result}, we will give a closed form expression for the sum 
\begin{eqnarray*}
\displaystyle\sum_{\lambda\in L_n}f_{\lambda}.
\end{eqnarray*}

Towards this goal, the following proposition is useful for our purposes. The claims therein can be proved easily using the hooklength formula. 
\begin{Proposition}\label{prop:use of hooklength formula}
Given $n\geq 2$, we have that $f_{(n,n)}=f_{(n,n-1)}=C_n$. Furthermore, 
\begin{eqnarray*}
&&f_{(n,n-1,1)}=\displaystyle\left( \frac{(n-1)(n+1)}{2n+1}\right)C_{n+1},\\
&&f_{(n-1,n-1,1)}=\displaystyle\left( \frac{n-1}{2}\right)C_{n}.
\end{eqnarray*}
\end{Proposition}
The results obtained about the Kronecker products  $s_{(n,n-1,1)}*s_{(n,n)}$ and $s_{(n-1,n-1,1)}*s_{(n,n-1)}$ in Sections \ref{section:first case} and \ref{section: second case} if interpreted in terms of characters of the symmetric group imply the following relations 
\begin{equation} \displaystyle\sum_{\substack{\lambda \vdash 2n\\l(\lambda) \leq 4}}(d_{\lambda}-1)f_{\lambda} + \displaystyle\sum_{\substack{\lambda \vdash 2n\\l(\lambda)=5 \\ \lambda_{5}=1}}f_{\lambda} = f_{(n,n-1,1)}f_{(n,n)}, \label{even2} \end{equation}
and
\begin{equation}\displaystyle\sum_{\substack{\lambda \vdash 2n-1\\l(\lambda) \leq 4}}(d_{\lambda}-1)f_{\lambda} + \displaystyle\sum_{\substack{\lambda \vdash 2n-1\\l(\lambda)=5 \\ \lambda_{5}=1}}f_{\lambda} = f_{(n-1,n-1,1)}f_{(n,n-1)}.\label{odd2}\end{equation}
Now, define $\sigma_{k}(n)$ as follows.
\begin{eqnarray*}\sigma_{k}(n)=\displaystyle\sum_{\lambda\vdash n, l(\lambda)\leq k} d_{\lambda}f_{\lambda}\end{eqnarray*}  Next, we will give a simple expression for $\sigma_{k}(n)$.
\begin{Theorem}\label{theorem:relating sigma and tau}
Given positive integers $m$ and $k$, we have
\begin{eqnarray*}\sigma_{k}(m) &=& \tau_{k}(m+1)-\tau_{k-1}(m) .\end{eqnarray*}
\begin{proof}
By definition we have that
\begin{eqnarray}\tau_{k}(m+1)&=& \displaystyle\sum_{\substack{\lambda \vdash m+1\\l(\lambda) \leq k}}f_{\lambda}.
\end{eqnarray}
Using \cite[Lemma 2.8.2]{Sa}, which says $f_{\lambda}=\sum_{\mu \prec \lambda}f_{\mu}$, we obtain the following sequence of equalities
\begin{eqnarray*} \displaystyle\sum_{\substack{\lambda \vdash m+1\\l(\lambda) \leq k}}f_{\lambda} &=& \displaystyle\sum_{\substack{\lambda \vdash m+1\\l(\lambda) \leq k}}\sum_{\mu \prec \lambda}f_{\mu}\nonumber\\ &=&\displaystyle\sum_{\substack{\mu \vdash m\\l(\mu) \leq k}}\sum_{\substack{\lambda \succ \mu \\ l(\lambda) \leq k}}f_{\mu} \nonumber \\ &=& \displaystyle\sum_{\substack{\mu \vdash m\\l(\mu) \leq k-1}}(d_{\mu}+1)f_{\mu}+\sum_{\substack{\mu \vdash m \\ l(\mu) = k}}d_{\mu}f_{\mu}\nonumber \\ &=& \displaystyle\sum_{\substack{\lambda \vdash m \nonumber\\l(\lambda) \leq k}}d_{\lambda}f_{\lambda}+\tau_{k-1}(m)\\ &=& \sigma_{k}(m)+\tau_{k-1}(m) .\end{eqnarray*}
Thus the claim is established.\end{proof}
\label{thesigmatheorem}
\end{Theorem}
Using Theorem \ref{theorem:relating sigma and tau} in conjunction with known expressions for $\tau_{3}(n)$ and $\tau_{4}(n)$ given in \eqref{eqn:motzkin} and \eqref{eqn:catalan} respectively, we get the following corollary.
\begin{Corollary}\label{thesigmafour}
Given a positive integer $n$, we have
\begin{eqnarray*}
\sigma_{4}(n)=C_{\lfloor\frac{n}{2}\rfloor +1}C_{\lceil\frac{n}{2}\rceil +1}-M_{n}.
\end{eqnarray*}
\end{Corollary}
Now we come to our main enumerative result which makes use of \eqref{even2} and \eqref{odd2}. It relates to a specific case of counting standard Young tableau with a fixed height.
\begin{Theorem}\label{theorem:main combinatorial result}
Given a positive integer $k\geq 3$, we have
\begin{eqnarray*}
\displaystyle\sum_{\lambda\in R_k}f_{\lambda} = \displaystyle\frac {\lfloor \frac{k+1}{2} \rfloor (\lceil \frac{k+1}{2} \rceil +1)}{k+1}C_{\lfloor \frac{k+1}{2} \rfloor}C_{\lceil \frac{k+1}{2} \rceil} - C_{\lfloor \frac{k}{2} \rfloor +1}C_{\lceil \frac{k}{2} \rceil +1}+M_{k}.\nonumber\\
\end{eqnarray*}
\begin{proof}
We will treat the cases where $k$ is odd and $k$ is even separately. Firstly assume $k=2n$ for some integer $n \geq 2$. Then \eqref{even2} implies
\begin{eqnarray}
\displaystyle\sum_{\lambda\in L_{2n}}f_{\lambda}&=& f_{(n,n-1,1)}f_{(n,n)}-\displaystyle\sum_{\substack{\lambda \vdash 2n\\l(\lambda) \leq 4}}(d_{\lambda}-1)f_{\lambda} \nonumber \\ &=& f_{(n,n-1,1)}f_{(n,n)}-\sigma_{4}(2n)+\tau_{4}(2n).
\end{eqnarray}
Using Corollary \ref{thesigmafour} for $\sigma_{4}(2n)$, \eqref{eqn:catalan} for $\tau_{4}(2n)$ and Proposition \ref{prop:use of hooklength formula} for $f_{(n,n-1,1)}$ and $f_{(n,n)}$ in the right hand side of the above equation, we obtain
\begin{eqnarray}
\displaystyle\sum_{\lambda\in L_{2n}}f_{\lambda} &=& \left(\frac{n^2-1}{2n+1}\right)C_{n}C_{n+1}-C_{n+1}^{2}+M_{2n}+C_{n}C_{n+1}\nonumber \\ &=& \left(\frac{n(n+2)}{2n+1}\right)C_{n}C_{n+1}-C_{n+1}^{2}+M_{2n}.
\end{eqnarray}
Now, assume $k=2n-1$ for $n \geq 2$. Then \eqref{odd2} implies
\begin{eqnarray}
\displaystyle\sum_{\lambda \in L_{2n-1}}f_{\lambda}&=& f_{(n-1,n-1,1)}f_{(n,n-1)}-\displaystyle\sum_{\substack{\lambda \vdash 2n-1\\l(\lambda) \leq 4}}(d_{\lambda}-1)f_{\lambda} \nonumber \\ &=& f_{(n-1,n-1,1)}f_{(n,n-1)}-\sigma_{4}(2n-1)+\tau_{4}(2n-1).
\end{eqnarray}
Using Corollary \ref{thesigmafour} for $\sigma_{4}(2n-1)$, \eqref{eqn:catalan} for $\tau_{4}(2n-1)$ and Proposition \ref{prop:use of hooklength formula} for $f_{(n-1,n-1,1)}$ and $f_{(n,n-1)}$ in the right hand side of the above equation, we obtain
\begin{eqnarray}
\displaystyle\sum_{\lambda\in L_{2n-1}}f_{\lambda} &=& \left(\frac{n-1}{2}\right)C_{n}^{2}-C_{n}C_{n+1}+M_{2n-1}+C_{n}^2 \nonumber \\ &=& \left(\frac{n+1}{2}\right)C_{n}^{2}-C_{n}C_{n+1}+M_{2n-1}.
\end{eqnarray}
The claim is now a unified way of rewriting the formulae obtained in the two cases, $k=2n$ and $k=2n-1$.
\end{proof}
\end{Theorem}

\section{Further directions}\label{section:future avenues}
It is clear from the techniques we have used here that if we know a combinatorial rule for computing the Kronecker product of Schur functions indexed by partitions of length at most $k$, for some positive integer $k$, then we can give a description of the Kronecker product of Schur functions indexed by partitions of length $k+1$, where the smallest part is $1$ or $2$. For example, if we wish to compute $s_{(n+k,n-k-1,1)}*s_{(n,n)}$ for $n \geq k+2$ and $k \geq 1$, then performing calculations akin to those described earlier, we obtain
\begin{eqnarray}
\langle s_{(n+k,n-k-1,1)}*s_{(n,n)},s_{\theta}\rangle &=& \displaystyle\sum_{j=0}^{k} (-1)^{k+j}\langle s_{(1)}(s_{(n+j,n-j)}*s_{(n,n)}),s_{(1)}s_{\theta} \rangle \nonumber \\ &&-\langle (s_{(n+k+1,n-k-1)}+s_{(n+k,n-k)})*s_{(n,n)},s_{\theta}\rangle \label{conclusioneq}\end{eqnarray}
Note that the case $k=0$ has already been dealt with in Section \ref{section:first case}. There is a combinatorial rule for computing Kronecker products of the form $s_{(n+j,n-j)}*s_{(n,n)}$, as described in \cite{BvWZ}. Thus, in theory, we can compute $\langle s_{(n+k,n-k-1,1)}*s_{(n,n)},s_{\theta}\rangle$ given \eqref{conclusioneq}. 

Another area worth investigating is counting standard Young tableaux with added constraints as described below. Given $k\geq 0$ and $i,n\geq 1$, consider the set 
\begin{eqnarray*}
S(k,i,n)=\{\lambda \vdash n: \lambda_{k+1}=i \text{ and } l(\lambda)= k+1\}.
\end{eqnarray*}
Let  $\rho_{k,i}(n)$ be defined by
\begin{eqnarray}
\rho_{k,i}(n)=\displaystyle\sum_{\lambda \in S(k,i,n)} f_{\lambda}.
\end{eqnarray}
Note that what we have enumerated in Theorem \ref{theorem:main combinatorial result} is $\rho_{4,1}(n)$. 
The numbers $\rho_{k,i}(n)$ provide a refinement of the sequence $\tau_{k}(n)$  and their enumeration may yield new combinatorial identities.

\end{document}